\documentclass[a4paper,reqno]{amsart}
\usepackage{kpfonts,baskervald}

\usepackage[T1]{fontenc}
\usepackage[utf8x]{inputenc}
\usepackage[english]{babel}
\usepackage{graphicx}
\usepackage{amsmath,mathrsfs,  amssymb,slashed,url,bm,upgreek,amsthm,esint}
\usepackage{graphicx,enumerate,comment}
\usepackage[hypertexnames=false]{hyperref}
\newcounter{intro}

\newtheorem{theo}[intro]{Theorem}

\newtheorem{thm}{Theorem}[section]
\newtheorem{lem}[thm]{Lemma}
\newtheorem{prop}[thm]{Proposition}
\newtheorem{cor}[thm]{Corollary}

\newtheorem{rem}[thm]{Remark}

\newtheorem*{merci}{Acknowledgements}
\usepackage{xcolor}
\DeclareTextFontCommand{\dc}{\color{purple}}

\newcommand{\cref}[1]{Corollary~\ref{#1}}

\newcommand{\pref}[1]{Proposition~\ref{#1}}
\newcommand{\rref}[1]{Remark~\ref{#1}}
\newcommand{\tref}[1]{Theorem~\ref{#1}}

\def\tK{\text{K}}
\def\mcC{\mathscr C}

\def\cC{\mathcal C}

\def\cH{\mathcal H}

\def\cP{\mathcal P}
\def\cS{\mathcal S}

\def\N{\mathbb N}

\def\bD{\mathbb D}

\def\Z{\mathbb Z}
\def\R{\mathbb R}

\def\vol{\mathrm{v}}

\DeclareMathOperator{\dist}{d}

\DeclareMathOperator{\Geo}{Geo}
\DeclareMathOperator{\CD}{CD}
\DeclareMathOperator{\RCD}{RCD}
\DeclareMathOperator{\MCP}{MCP}

\DeclareMathOperator{\ricci}{Ricci}

\DeclareMathOperator{\dv}{dv}

\DeclareMathOperator{\argsinh}{argsinh}
\newcommand{\meas}{\mathfrak{m}}

\newcommand{\eps}{\ensuremath{\varepsilon}}
\newcommand{\Ric}{\ensuremath{\mbox{Ric}}}

\newcommand{\Ricm}{\ensuremath{\mbox{Ric}_{-}}}
\newcommand{\measrestr}{%
  \,\raisebox{-.127ex}{\reflectbox{\rotatebox[origin=br]{-90}{$\lnot$}}}\,%
  }

\newcommand{\weakto}{\rightharpoonup}

\newcommand{\di}{\mathop{}\!\mathrm{d}}

\def\cS{\mathcal S}

\title[Strong Kato limits can be branching]{Strong Kato limits can be branching}

\author{Gilles Carron}
\address{G. Carron, Nantes Université, CNRS, Laboratoire de Mathématiques Jean Leray, LMJL, UMR 6629, F-44000 Nantes, France.} 
\email{Gilles.Carron@univ-nantes.fr}

\author{Ilaria Mondello}
\address{I. Mondello, Université Paris Est Cr\'eteil, Laboratoire d'Analyse et Math\'ematiques appliqu\'es, UMR CNRS 8050, F-94010 Créteil, France, Institut Universitaire de France.}
\email{ilaria.mondello@u-pec.fr}

\author{David Tewodrose}
\address{D. Tewodrose, Vrije Universiteit Brussel, Department of Mathematics and Data Science, Pleinlaan 2, B-1050 Brussel, Belgium.}
\email{david.tewodrose@vub.be}

\begin{document}
\maketitle\begin{abstract}We provide an example of a non-collapsed strong Kato limit that is branching, essentially branching, and satisfies neither the $\CD(K,\infty)$ nor the $\MCP(K,N)$ conditions for any $K \in \mathbb{R}$ and $N \in [1,+\infty)$.  In particular, this space is not a Ricci limit space.  We also construct a compact non-collapsed strong Kato limit that cannot be obtained as Gromov--Hausdorff limit of closed Riemannian surfaces satisfying a uniform small $L^p$ bound à la Petersen--Wei.



\end{abstract}

\section{Introduction}
This paper aims to construct metric spaces arising as Gromov–Hausdorff limits of surfaces satisfying a strong Kato bound. In particular, we give two explicit examples which show substantial differences from the Gromov–Hausdorff limits and singular spaces with curvature–dimension bounds previously studied in the literature.

We briefly recall some notation and the framework in which a structure theory for Kato limits has been developed.  For a complete smooth Riemannian manifold $(M^n,g)$, we let $\Ricm\colon M\rightarrow [0,+\infty)$  be the negative part of the smallest eigenvalue of the Ricci tensor, so that $\ricci\ge -\Ricm\, g.$ The Kato constant of $(M,g)$ at time $T>0$ is defined as
$$\text{k}_T(M,g)=\sup_{x\in M} \int_0^T\left[ e^{-t\Delta_g}\Ricm\right](x)\di t$$
where $(e^{-t\Delta_g})_{t>0}$ is the heat semigroup generated by the Laplace--Beltrami operator $\Delta_g$. Since this semigroup preserves $L^\infty$ bounds,  any uniform Ricci curvature lower bound implies a linear control on the Kato constant:
\begin{equation}\label{eq:Ricci}
\ricci\ge -K\, g \qquad \Rightarrow \qquad  \forall T>0 : \, \, \, \text{k}_T(M,g)\le KT.
\end{equation}

By analogy with the study of Schrödinger operators in the Euclidean space, we say that a sequence of complete smooth Riemannian manifolds  $\{(M_\ell^n,g_\ell)\}_\ell$ satisfies a uniform Kato bound if there exist $T>0$ and \mbox{$f\colon (0,T]\rightarrow \R_+$} non-decreasing such that 
\begin{equation*}
f(t) \to 0 \mbox{ as } t \downarrow 0, \quad  \mbox{and} \quad \text{k}_t(M_\ell,g_\ell)\le f(t) \mbox{ for all } \ell \mbox{ and } t \in (0,T].
\end{equation*}
The bound is called \emph{strong} if in addition
$$\int_0^{T}\frac{f(t)}{t}dt<\infty.$$ 

A strong Kato bound is weaker than previous curvature assumptions considered in the study of Gromov--Hausdorff limits. Indeed, by \eqref{eq:Ricci},  it is implied by a uniform Ricci curvature lower bound. It also follows from a uniform smallness assumption on the $L^p$ norm of $\Ricm$ for some $p>n/2$, as introduced by P.~Petersen and G.~Wei; see \cite{RoseStollmann,CarronRose} or Lemma \ref{lem:LpSK}.  However, Kato potentials that do not lie in any such $L^p$ exist in $\mathbb{R}^n$ \cite{Bsimon}.

In \cite{CMT,CMT2, CMT23}, we described  analytic and geometric properties of Gromov--Hausdorff limits of sequences of Riemannian manifolds with a uniform Kato bound.  The richest structure theory applies to the class of non-collapsed strong Kato limits, namely pointed metric spaces $(X,\dist,o)$ that arise as pointed Gromov-Hausdorff limits of sequences $\{(M_\ell^n,g_\ell,o_\ell)\}_\ell$ satisfying a uniform strong Kato bound and a non-collapsing hypothesis, that is, there exists $v>0$ such that for all $\ell$,
$$\mathrm{v}_{\ell}\left(B_{\ell}\left(o_\ell,\sqrt{T}\right)\right)\ge v\, T^{\frac n2}.$$

Our structure results recovered those of J.~Cheeger and T.H.~Colding for Ricci limits \cite{ChCo97,ChCo00}, i.e.~limits of pointed Riemannian manifolds satisfying a uniform Ricci curvature lower bound. They also extended those obtained by P.~Petersen and G.~Wei \cite{PW1, PW2}, which apply to sequences with a uniform $L^p$ smallness assumption on the negative part of the Ricci curvature. In particular we proved that for any non-collapsed strong Kato limit $(X,\dist)$ the volume measures $\mathrm{v}_{\ell}$ converge to the $n$-dimensional Hausdorff measure $\mathcal{H}^n$. In \cite{CMT23}, we additionally established that the resulting metric measure space $(X,\dist,\cH^n)$ is bi-Lipschitz  to an $\RCD(K,N)$ space for finite $N$. 

This naturally lead to the following questions. First, which are the geometric differences between Kato limits and either Ricci limits or finite dimensional $\RCD$ spaces? Secondly, in analogy with the existence of Kato potentials that do not belong to $L^p$ for $p$ larger than half of the dimension: are there non-collapsed strong Kato limits that cannot arise as limits of manifolds with an $L^p$ bound?

\hfill

In this article, we exhibit examples of Kato limits which show that the class of non-collapsed strong Kato limits is strictly larger than both the one of finite-dimensional $\RCD$ spaces and the one covered by Petersen and Wei's theory.

Our first example is a non-collapsed strong Kato limit which is \emph{branching}, unlike finite dimensional $\RCD$ spaces.  Q.~Deng established in \cite{Deng} that these spaces are non-branching, building upon the work of T.H.~Colding and A.~Naber \cite{Colding_2012,Colding_2013} for limit geodesics in Ricci limits.  The example we construct below is even essentially branching and degenerate in the sense of \cite{Kell}, and therefore satisfies neither any $\CD(K,\infty)$ nor any finite-dimensional $\MCP(K,N)$ condition. 

Consider two Euclidean copies of $\R^2\setminus \bD^2$ glued along their boundaries, where $\bD^2\subset \R^2$ is the open unit disk. This produces the cylinder $$\mathcal{C} = \R \times (\R/(2 \pi \Z))$$ endowed with the singular warped Riemannian metric $$\widehat g=(dt)^2+(1+|t|)^2\,(d\theta)^2.$$  Let $\widehat \dist$ and $\widehat{\cH}^2$ be the associated distance and Hausdorff measure.  The following holds.

\begin{theo} 
\label{thm:exA}
The surface $(\mathcal{C},\widehat{d} )$ is a branching non-collapsed strong Kato limit. Moreover,  the metric measure space $(\mathcal{C},\widehat \dist,\widehat{\cH}^2)$ satisfies neither the $\CD(K,\infty)$ nor  the $\MCP(K,N)$ conditions for any $K \in \R$ and finite $N >1$.
\end{theo}

Note that $\widehat{g}$ is bi-Lipschitz equivalent to a smooth Riemannian metric with bounded curvature. More generally,  by \cite{CMT23},  any Gromov--Hausdorff limit of Riemannian surfaces with uniformly bounded Kato constant at a fixed time $T$ is bi-Lipschitz to an Alexandrov space with curvature bounded below, and such spaces are non-branching.

The key ingredient to establish that $(\mathcal{C},\widehat{d})$ is a non-collapsed strong Kato limit is a general criterion for Riemannian surfaces with low regularity (Theorem \ref{th:IBC}) for being such a limit space. Specifically,  if $(M,g_0)$ is a complete Riemannian surface with bounded Gauss curvature and $g_f=e^{2f}g_0$ is a conformal deformation where  $f$ is bounded, Lipschitz,  and has distributional Laplacian given by a Radon measure bounded from above,  then $(M,g_f)$ is a non-collapsed strong Kato limit.  Such surfaces are locally BIC (for Bounded Integral Curvature),  see \cite{Fillastre,Troyanov}.  This gives a broad class of low-regularity surfaces that are non-collapsed strong Kato limits,  opening the way to a wide range of new examples.

Our second example answers positively the previous question concerning $L^p$ bounds:  there exist non-collapsed strong Kato limits which cannot be obtained as limits of manifolds satisfying a uniform $L^p$ smallness condition à la Petersen-Wei. We construct such an example by performing an analogue of the previous gluing construction in a compact setting. Let two copies of $\mathbb{S}^2 \backslash \mathbb{D}_\alpha$ be glued together along their boundaries, where $ \mathbb{D}_\alpha$ is a geodesic disk of radius $\alpha \in (0,\pi/2)$ for the round metric. This yields a singular Riemannian metric $\overline{g}$ on $\mathbb{S}^2$. We show the following result, where $\epsilon(2,p)$ and $\mathbb{V}_{2,\kappa}$ are defined in Section \ref{sec:B}. For any function $f$ we let $f_-=-\min(0,f)$ be its negative part.

\begin{theo}
\label{prop:NotLp}
The surface  $(\mathbb{S}^2,\overline{g})$ is a non-collapsed strong Kato limit. Moreover,  it cannot arise as Gromov-Hausdorff limit of closed surfaces $\{(M_\ell, g_\ell)\}$ of diameters at most $D>0$ such that, for some $p>1$ and $\kappa \geq 0$,
$$\sup_\ell D^2\left( \fint_{M_\ell} (K_{g_\ell} +\kappa^2 )_-^p  \di\vol_{g_\ell} \right)^{1/p} \le \epsilon(2,p) \,  \left( \frac{D^{2}}{\mathbb{V}_{2,\kappa}(D)}\right)^{1/p} \, \cdot $$
\end{theo}

The paper is structured as follows. Section \ref{sec:criterion} establishes Theorem \ref{th:IBC}, the general criterion yielding a large class of two-dimensional non-collapsed strong Kato limits, including $(\mathcal{C},\widehat{d})$ and $(\mathbb{S}^2,\overline{g})$. Section \ref{sec:A} analyses the first example in detail, proving in particular its branching and non-$\CD$/non-$\MCP$ properties. Section \ref{sec:B} proves Theorem \ref{prop:NotLp} by means, among others, of the pre-compactness result of \cite{Debin} for compact surfaces: an analogous result in the complete case would similarly imply that $(\mathcal{C},\widehat{d})$ cannot be a limit of surfaces satisfying a uniform smallness $L^p$ bound.

\begin{merci}
The third author has been supported by Laboratoire de Mathématiques Jean Leray via the project Centre Henri Lebesgue ANR-11-LABX-0020-01,  by Fédération de Recherche Mathématiques de Pays de Loire via the project Ambition Lebesgue Loire,  and by the Research Foundation – Flanders (FWO) via the Odysseus II programme no.~G0DBZ23N. The second author was supported by the National Science Foundation under Grant No. DMS-1928930, while she was in residence at the Mathematical Sciences Research Institute in Berkeley, California, during the fall semester of 2024. The authors thank Andrea Mondino and Shouhei Honda for their stimulating questions raised during the School and Conference on Metric Measure Spaces, Ricci Curvature, and Optimal Transport, at Villa Monastero, Lake Como. They also thank Mauricio Che for pointing out the reference \cite{Kell}.
\end{merci}

\section{Conformal metric with low regularity and strong Kato limits}\label{sec:criterion}

In this section, $c(a_1,a_2, \ldots)$ denotes a generic positive constant depending on parameters $a_1,a_2, \ldots$ only. The value of this constant may change from line to line.

\subsection{The Kato constant of a smooth conformal metric}

The next key statement shows that, on a smooth surface with lower bounded curvature,  the Kato constant of a conformal metric is controlled upon mild assumptions on the conformal factor.

\begin{prop}\label{prop:Katoest} Let $(\Sigma,g_0)$ be a complete Riemannian surface whose Gaussian curvature is uniformly bounded from below: there exists $\kappa\ge 0$ such that
$$\tK_{g_0}\ge -\kappa^2.$$ 
Assume that $f\in \mathscr{C}^\infty(\Sigma)$ is a bounded function uniformly Lipschitz with respect to $g_0$, i.e.  there exist $L,\lambda >0$ such that
\begin{equation}\label{eq:bounded}
\|f\|_{L^\infty}\le L,
\end{equation}
\begin{equation}\label{eq:Lipschitz}
\forall x,y\in \Sigma\colon |f(x)-f(y)|\le \lambda d_{g_0}(x,y),
\end{equation}
and that the Gaussian curvature of $g_{f}=e^{2f}g_0$ is uniformly 
bounded from above: there exists $K\in \mathbb{R}$ such that
$$\tK_{g_f}\le K.$$
Then for any $T\le \kappa^{-2}$, the Kato constant of the metric $g_f$ satisfies 
$$\text{k}_T(\Sigma,g_f)\le c(\kappa,L,\lambda,K_+) \sqrt{T}.$$
\end{prop}

To prove this result, we need some estimates on the heat kernels of $g_0$ and $g_f$. Firstly,  since the metric $g_0$ has lower bounded curvature, the Li-Yau estimates provide (\cite{LiYau}):
\begin{equation}\label{E1} 
\forall x,y\in \Sigma, \forall t\in \left(0,\kappa^{-2}\right)\colon H_{g_0}(t,x,y)\le \frac{c(\kappa)}{\vol_{g_0} \left(B_{g_0}(x,\sqrt{t})\right) }\, e^{-\frac{d^2_{g_0}(x,y)}{5t}}.\end{equation}
Moreover,  the function $f$ is bounded \eqref{eq:bounded} so $g$ and $g_0$ are bi-Lipschitz,  hence similar estimates hold for the heat kernel of $(\Sigma,g_f)$ (\cite{Grigoryan92,Saloff_Coste_1992, Saloff-Coste}):
 \begin{equation}\label{E2} 
\forall x,y\in \Sigma, \forall t\in \left(0,\kappa^{-2}\right)\colon H_{g_f}(t,x,y)\le \frac{c(\kappa,L)}{\vol_{g_f} \left(B_{g_f}(x,\sqrt{t})\right) }\, e^{-\frac{d^2_{g_f}(x,y)}{5t}}.\end{equation}
We also know (see \cite[Lemma 3.6]{Hebisch_2001}) that for any $x\in \Sigma$ and any $t,R>0$ such that $0\le t\le \kappa^{-2}$ and $R>\sqrt{t}$:
\begin{equation}\label{IEext} 
\int_{\Sigma \setminus B_{g_f}(x,R)} H_{g_f}(t,x,y)\dv_{g_f}(y)\le c(\kappa,L)\ e^{-\tilde{c} \frac{R^2}{t}}\end{equation}where $\tilde{c}$ depends only on $\kappa$ and $L$. Lastly, according to the Li-Yau gradient estimate \cite{LiYau},  since $f$ is uniformly Lipschitz \eqref{eq:Lipschitz}, then $e^{-t\Delta_{g_0}}f$ is also uniformly Lipschitz:
\begin{equation}\label{liplip} 
\left\|\, de^{-t\Delta_{g_0}}f\right\|_{L^\infty}\le e^{\kappa^2 t}\lambda.
\end{equation}

We are now in a position to prove \pref{prop:Katoest}.

\begin{proof}
The transformation rule for the Gaussian curvature under a conformal change is 
\begin{equation}
\label{eq:CT}
\tK_{g_f}=e^{-2f}\tK_{g_0}+e^{-2f}\Delta_{g_0} f=e^{-2f}\tK_{g_0}+\Delta_{g_f} f.
\end{equation}
The hypothesis $\tK_{g_f}\le K$ implies that $\left(\tK_{g_f}\right)_+ \le K_+$ so that
$$\left(\tK_{g_f}\right)_- = \left(\tK_{g_f}\right)_+ - \tK_{g_f} \le K_+-\tK_{g_f}=K_+-e^{-2f}\tK_{g_0}-\Delta_{g_f} f.$$
Hence for $x\in \Sigma$ and $T\le \kappa^{-2}$, one gets
\begin{align*}
\int_{[0,T]\times \Sigma}\!\!\!  H_{g_f}(t,x,y)\left(\tK_{g_f}\right)_-(y)\dv_{g_f}(y)\di t\le &T\left(K_++\kappa^2 e^{2L}\right)\\
&-\int_{[0,T]\times \Sigma}\!\! \!\!H_{g_f}(t,x,y)\Delta_{g_f} f(y)\dv_{g_f}(y) \di t.\end{align*}
But
\begin{align*}-\int_{[0,T]\times \Sigma}\!\! \!\!H_{g_f}(t,x,y)\Delta_{g_f} f(y)\dv_{g_f}(y) \di t&=\int_0^T \frac{\partial}{\partial t} \left(e^{-t\Delta_{g_f}}f\right)(x)  \di t\\
&=\left(e^{-T\Delta_{g_f}}f\right)(x)-f(x).\end{align*}
As the surface $(\Sigma,g_f)$ is stochastically complete \cite[Theorem 2.1]{Hebisch_2001}, we have that 
$$\int_\Sigma H_{g_f}(T,x,y) \dv_{g_f}(y)=1,$$ so that 
$$\left|\left(e^{-T\Delta_{g_f}} f \right)(x)-f(x)\right|\le \int_\Sigma H_{g_f}(T,x,y)\left|f(y)-f(x)\right|\dv_{g_f}(y).$$
By assumption
$$\left|f(y)-f(x)\right|\le e^{L} \lambda d_{g_f}(x,y)$$ and we are left with the estimate of 
$$\int_\Sigma H_{g_f}(T,x,y)\, d_{g_f}(x,y)\dv_{g_f}(y).$$
We have
\begin{align*}\int_\Sigma H_{g_f}(T,x,y)\,  d_{g_f}(x,y)\dv_{g_f}(y)&\le\int_{B_{g_f}(x,\sqrt{T})} H_{g_f}(T,x,y)\,  d_{g_f}(x,y)\dv_{g_f}(y) \\
&\ \ \ +\int_{\Sigma\setminus B_{g_f}(x,\sqrt{T})} H_{g_f}(T,x,y)\,  d_{g_f}(x,y)\dv_{g_f}(y)\\
&\le \sqrt{T}+\int_{\Sigma\setminus B_{g_f}(x,\sqrt{T})} H_{g_f}(T,x,y)\,  d_{g_f}(x,y)\dv_{g_f}(y).
\end{align*}
Let us introduce the  Stieltjes measure associated to the non-decreasing function
$$I(r)=-\int_{\Sigma\setminus B_{g_f}(x,r)}H_{g_f}(T,x,y)\dv_{g_f}(y).$$
Using \eqref{IEext},  we get that $-I(r)\le  c(\kappa, L)e^{-c\frac{r^2}{T}}$ for any $r\ge \sqrt{T}$, hence we can justify the integration by part formula and the following estimates:
\begin{align*}
\int_{\Sigma\setminus B_{g_f}(x,\sqrt{T})} H_{g_f}(T,x,y)\,  d_{g_f}(x,y)\dv_{g_f}(y)&=\int_{\sqrt{T}}^\infty r \di I(r)\\
&= -\sqrt{T}I(\sqrt{T})-\int_{\sqrt{T}}^{+\infty}I(r)\di r\\
&\le \sqrt{T}+c(\kappa, L)\int_{\sqrt{T}}^{+\infty}e^{-c\frac{r^2}{T}}\,\di r\\
&\le\sqrt{T}\left(1+c(\kappa, L)\int_{0}^{+\infty}e^{-c\rho^2}\di\rho\right),
\end{align*}
so that we obtain
\begin{equation}\label{estfinal}\int_\Sigma H_{g_f}(T,x,y)\,  d_{g_f}(x,y)\dv_{g_f}(y)\le c(\kappa,L)\sqrt{T}.\end{equation}
We eventually get the estimate
$$\text{k}_T(\Sigma, g_f)\le c(\kappa,L)\lambda \sqrt{T}+\left(K_++\kappa^2 e^{2L}\right) T.$$
\end{proof}

\subsection{Conformal metric with low regularity}
Building upon Proposition \ref{prop:Katoest}, we prove the following.

\begin{thm}\label{th:IBC}Let $(\Sigma,g_0)$ be a complete smooth Riemannian surface whose Gaussian curvature is uniformly bounded: there exists $\kappa \ge 0$ such that
$$\left|\tK_{g_0}\right|\le \kappa^2.$$ 
Assume that $f\in \mathscr{C}^0(\Sigma)$ is a bounded function  uniformly Lipschitz with respect to $g_0$, i.e.~there exist $L,\lambda \ge 0$ such that \eqref{eq:bounded} and \eqref{eq:Lipschitz} hold, and assume that the distributional Laplacian $\Delta_{g_0} f$ is a Radon measure uniformly bounded from above: there exists $M \in \mathbb{R}$ such that
\begin{equation}\label{eq:boundedLaplacian}
\Delta_{g_0} f\le M \text{ weakly.}\end{equation}
Then $(\Sigma,g_f)$ is a non-collapsed strong Kato limit.
\end{thm}
\begin{rem}\label{rem:Bac}  Since the Gaussian curvature of $g_0$ is uniformly bounded, the hypothesis on the Laplacian of $f$ is equivalent to demanding that the Gaussian curvature $\tK_{g_f}$ is a Radon measure --- so that the metric $g_f$ is locally BIC (\cite{Fillastre,Troyanov}) --- and that
$$\tK_{g_f}\le M \text{ weakly}.$$
\end{rem}

Let us prove Theorem \ref{th:IBC}.

\proof We are going to find $\{f_\epsilon\}\subset\mathscr{C}^\infty(\Sigma)$  such that 
\begin{enumerate}[i)]
\item $f_\epsilon\to f\text{ uniformly  as $\epsilon \downarrow 0$,}$
\item $f_\epsilon$ is uniformly bounded and Lipschitz, i.e.~there is $P>0$ independent of $\epsilon$ such that
$$\|f_\epsilon\|_{L^\infty}\le P\text{ and } \|df_\epsilon\|_{L^\infty}\le P ;$$
\item  the Gaussian curvature of $g_{\epsilon}=e^{2f_\epsilon}g_0$ is uniformly 
bounded from above: there is some constant $C$ independent of $\epsilon,$ such that
$$\tK_{g_\epsilon}\le C.$$
\end{enumerate}
The first property implies that for any $o\in \Sigma$:
$$\left(\Sigma, d_{g_\epsilon},o\right)\stackrel{\text{pGH}}{\longrightarrow}\left(\Sigma, d_{g},o\right)\text{ as }\epsilon\to 0.$$
The second one also implies that we have a non-collapsing hypothesis on all the metrics $g_\epsilon$ and the second and the third properties yield a uniform control of the Kato constant of $(\Sigma, g_\epsilon)$ thanks to \pref{prop:Katoest}.\\
We define $f_\epsilon$ to be the heat regularization of $f$ at time $\epsilon$:
$$f_\epsilon=e^{-\epsilon\Delta_{g_0}}f.$$
Since $e^{-\epsilon\Delta_{g_0}}$ preserves $L^\infty$ bounds, we get from \eqref{eq:bounded} that 
$$\|f_\epsilon\|_{L^\infty}\le L.$$ Moreover, thanks to \eqref{liplip}, we have
$$\|\,df_\epsilon\|_{L^\infty}\le e^{\kappa^2 \epsilon} \lambda.$$
Thus ii) holds.  For any $x\in \Sigma,$
$$\left| f_\epsilon(x)-f(x)\right|\le \int_\Sigma H_{g_0}(\epsilon,x,y)\left|f(y)-f(x)\right|\dv_{g_0}(y)\le \lambda\int_\Sigma H_{g_0}(\epsilon,x,y)\, d_{g_0}(x,y)\dv_{g_0}(y),$$
and the proof of \pref{prop:Katoest} (see \eqref{estfinal}) shows that when $\epsilon\le \kappa^{-2}$ then
$$\int_\Sigma H_{g_0}(\epsilon,x,y)\, d_{g_0}(x,y)\dv_{g_0}(y)\le c(\kappa, L) \sqrt{\epsilon}.$$
Therefore, when $\epsilon$ tends to $0$, the sequence $\left(f_\epsilon\right)$ converges to $f$ uniformly, so i) is checked.
We need to ensure now that iii) holds, that is to say, that the curvature  of $g_{\epsilon}=e^{2f_\epsilon}g_0$ is uniformly 
bounded from above. Recall that
$$\tK_{g_\epsilon}=e^{-2f_\epsilon}\left(\tK_{g_0}+\Delta_{g_0}f_\epsilon\right).$$
But \eqref{eq:boundedLaplacian} yields that
$$\Delta_{g_0}f_\epsilon=e^{-\epsilon\Delta_{g_0}}\Delta_{ g_0} f\le M$$ so that 
$$\tK_{g_\epsilon}\le e^{2L}\left( K+M\right).\eqno \qedhere$$
\endproof

\begin{rem}
In \tref{th:IBC}, it is easy to verify that the same argument applies when $f$ is assumed to be uniformly Hölder rather than uniformly Lipschitz.
\end{rem}

\section{Proof of Theorem \ref{thm:exA}}\label{sec:A}

This section addresses Theorem \ref{thm:exA} that provides a branching non-collapsed strong Kato limit. We recall that for any function $f$ we denote by $f_-=-\min(0,f)$ its negative part and in the remainder of the paper we write $f_+=\max(0,f)$ for its positive part.

\subsection{The example is a non-collapsed strong Kato limit}

\begin{prop}
The space $(\cC,\widehat g)$ is a a non-collapsed strong Kato limit.
\end{prop}

\begin{proof}
The Gaussian curvature of a Riemannian metric on $\cC$ of the form
\[
g=(dr)^2+J^2(r)\,(d\theta)^2
\]
is given by:
\begin{equation}\label{eq:formula}
 K_g(r,\theta) = - \frac{J''(r)}{J(r)} \qquad \forall (r, \theta) \in \cC.
\end{equation}
The latter has to be intended in the weak sense when $J\in W_{\text{loc}}^{1,1}(\mathbb{R})$. Endow $\cC$ with the smooth Riemannian metric
$$g_0=(dr)^2+(1+r^2)\,(d\theta)^2.$$
Then \eqref{eq:formula} yields that for any $(r,\theta) \in \cC$,
\[
K_{g_0}(r,\theta) = -\frac{1}{(1+r^2)^2} \le 0.
\]
Moreover, it is easy to check that, up to reparametrization $t=r-1+\mathrm{sgn}(r)\sqrt{1+r^2}$ where $\mathrm{sgn}$ is the sign function, the metric $\widehat g$ is conformal to
$g_0$ with conformal factor 
$$f(r,\theta)= \argsinh(|r|)-\log\left(\sqrt{1+r^2}\right)=\log\left(1+\frac{|r|}{\sqrt{1+r^2}}\right).$$
The function $f$ is clearly bounded and uniformly Lipschitz.  Letting $\delta_0$ denote the Dirac measure on $\mathbb{R}$ centered at $0$,  we have
\[
- \frac{ \frac{d^2}{dr^2}(1+|r|)}{1+|r|} = - 2\delta_0 \quad \text{weakly.}
\]
Then \eqref{eq:formula} yields that for any $\varphi \in \cC_c(\cC)$, 
\[
\iint_{\R\times \left(\R/2\pi\Z\right)} \varphi(r,\theta) \di  K_{\widehat{g}}(r,\theta) =  - 2\int_{\R/2\pi\Z}\varphi(0,\theta) \di \theta,
\]
hence
\[
K_{\widehat g} = - 2 \di \mathcal{H}^1 \measrestr \cS 
\]
where $\di \mathcal{H}^1$ is  the  Euclidean one-dimensional Hausdorff measure on $\cC$,  $\measrestr$ means measure restriction, and $\cS = \{0\}\times \left(\R/2\pi\Z\right) \subset \cC$. In particular,   the distributional Gaussian curvature of $\widehat g$ is a bounded non-positive Radon measure. Using \rref{rem:Bac}, we can apply \tref{th:IBC} to conclude that $(\cC,\widehat g)$ is a non-collapsed strong Kato limit.
\end{proof}

\subsection{The example is branching}

If $(X,\dist)$ is a complete and separable metric space, a geodesic is a map $\upgamma\colon [0,1]\rightarrow X$ such that :
\begin{equation}\label{eq:defgeo}
\forall s,s'\in [0,1]\colon\ d(\upgamma(s),\upgamma(s'))=|s-s'|d(\upgamma(0),\upgamma(1)).
\end{equation} A geodesic $\gamma$ is called branching if there  exist another geodesic $\upbeta :  [0,1] \to X$ and a number $c\in (0,1)$ such that
$\upgamma(1)\not=\upbeta(1)$ but $\upgamma(t)=\upbeta(t)$ for any $t\in [0,c]$.  We say that $(X,d)$ is branching if it admits a branching geodesic. 

\begin{prop}
The space $(\mathcal{C},\widehat \dist )$ is branching. 
\end{prop}

\begin{proof}
Set $\Omega_\pm=\{(t,\theta)\in \mathcal{C} \colon \pm t\ge 0\}.$ Note that the maps 
\[
\begin{array}{ccccl}
\pi_\pm & : & \mathcal{C} & \to & \Omega_\pm\\
& &  (t,\theta) & \mapsto & (\pm t_\pm, \theta)
\end{array}
\]
are $1-$Lipschitz projections. Moreover, if $\gamma\colon [0,1]\rightarrow \mcC$ is a continuous rectifiable curve, then the length of $\pi_\pm(\gamma)$ is shorter that the  length of $\gamma$ with equality if and only if $\gamma([0,1])\subset\{t=0\}$. As a consequence,  a geodesic joining two points $p,q\in \Omega_+$ stays in $\Omega_+$.  Moreover, the maps
\[
\begin{array}{ccccl}
\Phi_\pm & :&  \Omega_\pm & \to & \R^2 \backslash \mathbb{D}^2 \\
& &  (t,\theta) & \mapsto & (1\pm t)(\cos \theta, \sin \theta)
\end{array}
\]
satisfy $\Phi_\pm^*g_e = \widehat{g}$, where $g_e$ is the Euclidean Riemannian metric.  Therefore, the geodesic joining $p,q\in \Omega_+$ is given by the pre-image by $\Phi_+$ of the segment $[\Phi_+(p),\Phi_+(q)] \subset \R^2$, provided this segment does not intersect with $\mathbb{D}^2$.  From the $1-$Lipschitz projection 
$$(t,\theta)\in\cC\mapsto (0,\theta)\in \partial\Omega_+,$$
we also get that if a geodesic $\gamma\colon [0,1]\rightarrow \Omega_+$ joining two points $p,q\in \Omega_+$ touches $\partial \Omega_+$ at two distinct points $$\gamma(a),\gamma(b)\in \partial \Omega_+, \qquad 0 \le a<b \le 1,$$ then it remains in $\partial \Omega_+$  between the two points :
$$\gamma([a,b])\subset \partial\Omega_+.$$
This implies that the unique geodesic  joining $(0,0)$ to  $(0,\pi)$ is $\gamma(s)=(0,s),\ s\in [0,\pi]$ and that $\gamma$ branches off at any interior point when one chooses to follow the half-line in $\Omega_+$ tangent to  $\partial\Omega_+$, see Figure 1.
\end{proof}

\begin{figure}
\includegraphics[height=3cm]{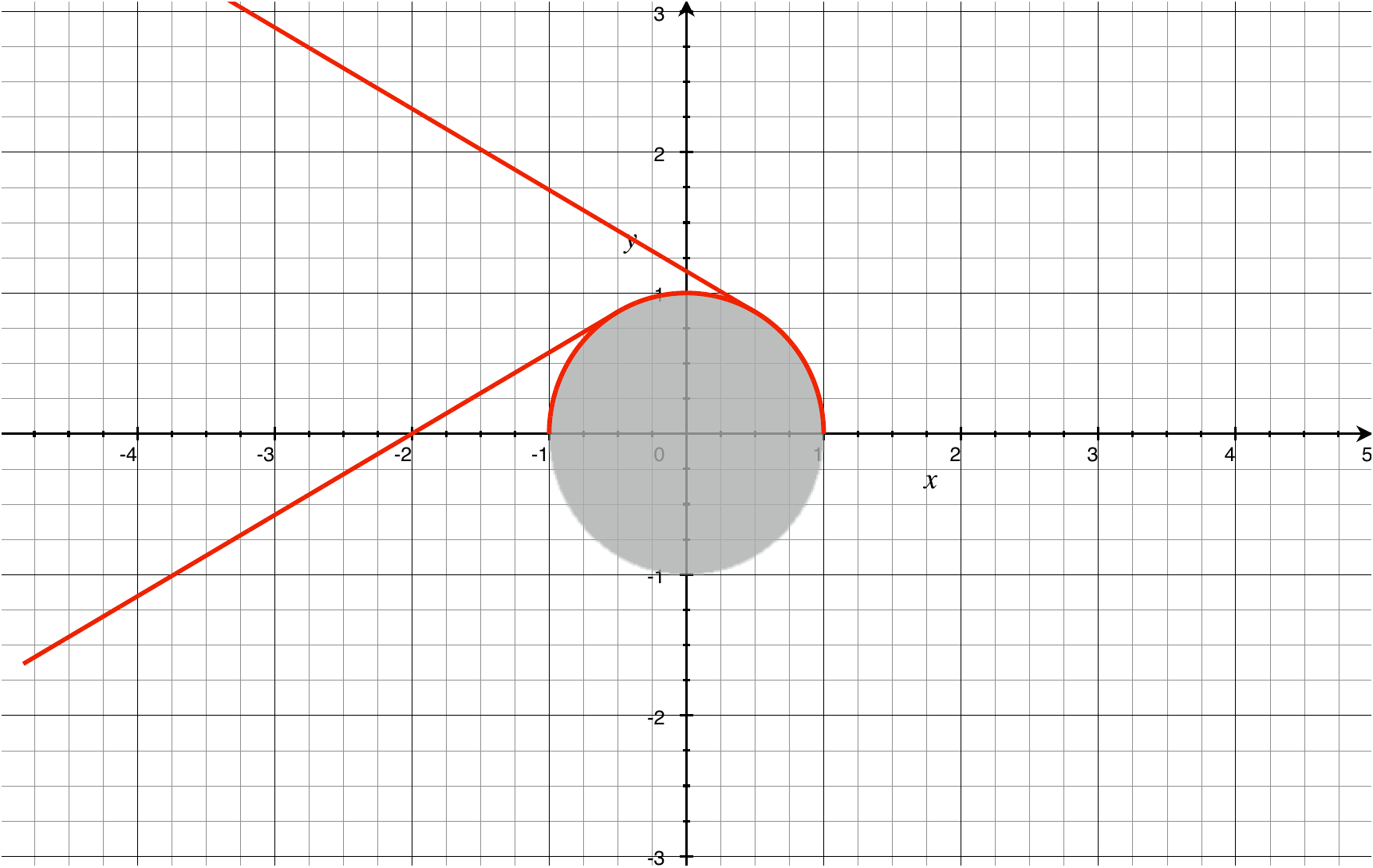}
\caption{Branching geodesics.}
\end{figure}

\begin{rem}\label{rem:symmetry} This surface has an isometric symmetry $p=(t,\theta)\mapsto S(p)=(-t,\theta)$. Hence   the minimizing geodesic $\gamma$ can also branch off from any of its interior points through a tangent  half-line in  $\Omega_-$. \end{rem}

\subsection{The example is essentially branching and not CD}

Let $(X,d)$ be a complete and separable metric space.  Denote by $\mathrm{Geo}(X)$ the set of geodesics in $X$ endowed with the uniform metric inherited from $ \cC([0,1],X)$.  We say that $X$ is geodesic if any pair of points $x,y \in X$ can be linked by a geodesic, that is to say,  there exists $\gamma \in \mathrm{Geo}(X)$ such that $\gamma(0)=x$ and $\gamma(1)=y$.  A set $\Gamma \subset \mathrm{Geo}(X)$ is called non-branching if for any $\gamma,\beta \in \Gamma$, if there exists $t \in (0,1)$ such that $\gamma(s)=\beta(s)$ for any $s \in(0,t)$, then $\gamma=\beta$.  We assume from now on that $(X,d)$ is geodesic. 

Denote by $\mathcal{P}(X)$ the set of all probability measures on $X$, and by $\mathcal{P}_2(X)$ the subset of those probability measures having a finite second moment. Recall that the $L^2$ Wasserstein distance between $\mu_0,\mu_1 \in \mathcal{P}_2(X)$ is defined as
\begin{equation}\label{eq:ot}
W_2(\mu_0,\mu_1) = \inf_{\sigma \in \Pi(\mu_0,\mu_1)} \left( \int_{X\times X} d^2(x,y) \di \sigma(x,y) \right)^{1/2}
\end{equation}
where $\Pi(\mu_0,\mu_1)$ is the set of probability measures on $X\times X$ with first and second marginals equal to $\mu_0$ and $\mu_1$ respectively. The so-called Wasserstein space $(\cP_2(X),W_2)$ is a complete, separable and geodesic metric space. More precisely,  for any $t \in [0,1]$, let $e_t  : \mathrm{Geo}(X) \to X$ be  mapping $\gamma$ to $\gamma(t)$. Any element $\pi \in \cP(\cC([0,1]),X)$ induces a curve $(\mu_t)_{t \in [0,1]}$ in $\cP_2(X)$ defined by $\mu_t = (e_t)_\# \pi$. This curve is a geodesic in $(\cP_2(X),W_2)$ if and only if $\pi(\Geo(X))=1$, in which case we say that $\pi$ is an optimal geodesic plan between $\mu_0$ and $\mu_1$. We write $\mathrm{OptGeo}(\mu_0,\mu_1)$ for the set of such plans.

\begin{figure}
\includegraphics[height=3.5cm]{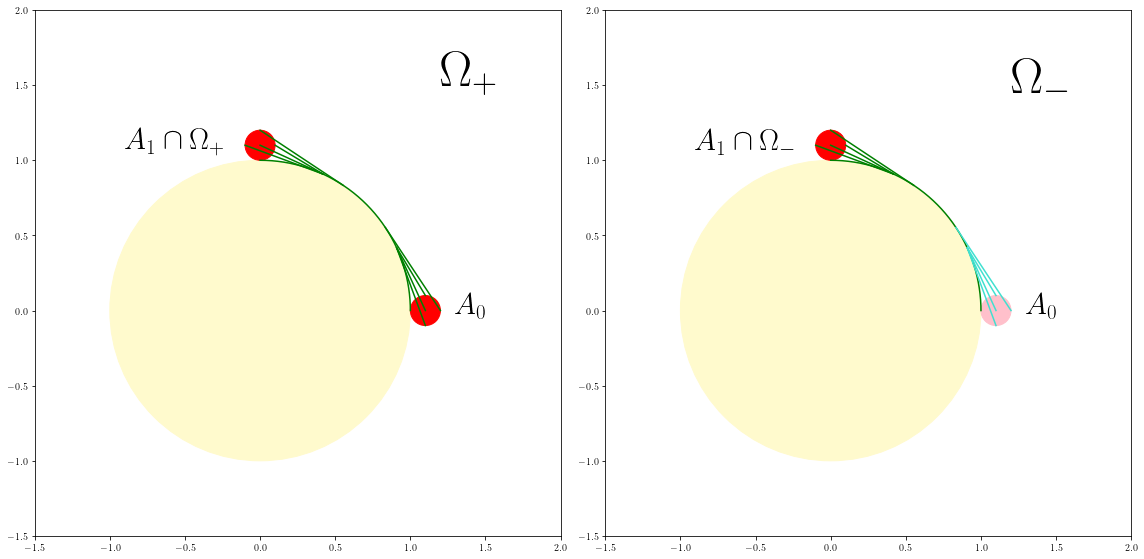}
\caption{Branching set of geodesics between $A_0$ and $A_1$.}
\end{figure}

Assume that $(X,d)$ supports a Radon measure $\meas$ which is finite on bounded sets.  Write $\cP_2^{\text{ac}}(X,\meas)$ for the subset of $\cP_2(X)$ whose elements are absolutely continuous with respect to $\meas$. Then $(X,d,\meas)$ is called essentially non-branching if for any couple $\mu_0, \mu_1 \in \mathcal{P}_2^{\text{ac}}(X,\meas)$, any $\pi \in \mathrm{OptGeo}(\mu_0,\mu_1)$ is concentrated in a set of non-branching geodesics.  This property was introduced in \cite{RajalaSturm} as a weaker measure-theoretical version of the non-branching property.  Note that there exist branching spaces which are essentially non-branching, see e.g.~\cite{OhtaExamples}. For brevity, we say that $(X,d,\meas)$ is essentially branching if it is not essentially non-branching.

\begin{prop}
The space $(\cC,\widehat g,\widehat \cH^2)$ is essentially branching. 
\end{prop}
\begin{proof}
For $\eps = 0.1$,  consider $o_0  = (\eps,0)  \in \Omega_+$ and $o_1^\pm = (\pm \eps, \pi/2)\in \Omega_\pm$. Set $A_0 =  \bar{B}_{\eps}(o_0)  \subset \Omega_+$ and $A_1=  \bar{B}_{\eps}(o_1^+) \cup  \bar{B}_{\eps}(o_1^-)$. For any $i \in \{0,1\}$, let $\mu_i$ be the restriction of $\widehat \cH^2$ to $A_i$ divided by $\widehat \cH^2(A_i)$. Then $\mu_0,\mu_1 \in \cP_2^{\text{ac}}(\cC,\widehat \cH^2)$.  Set $Y= \{(t,\theta) \in \mathcal{C} : -\pi/4 < \theta < 3\pi/4\}$ and note that for any $x,y \in Y$ there exists a unique geodesic $\gamma_{xy}$ such that $\gamma_{xy}(0)=x$ and $\gamma_{xy}(1)=y$, and $\gamma_{xy}([0,1]) \subset Y$ in addition. Consider the map $F : Y \times Y \to \Geo(\cC)$ defined by $F(x,y)=\gamma_{xy}$. Let $\sigma$ be an optimal plan between $\mu_0,\mu_1$, that is to say, a minimizer of \eqref{eq:ot}. Then the dynamic plan $\pi = F_\# \sigma$ belongs to $\mathrm{OptGeo}(\mu_0,\mu_1)$ and is concentrated in $\Gamma = \{\gamma \in \Geo(\cC) : \gamma(0) \in A_0 \text{ and } \gamma(1) \in A_1\}.$ But any $\gamma \in \Gamma$ is branching. Indeed,  set $x = \gamma(0)$, $y = \gamma(1)$, and $\tilde{y} = S(y)$ where $S$ is the isometric symmetry from Remark \ref{rem:symmetry}.  Let $\tilde{\gamma} \in \Gamma$ be the geodesic from $x$ to $\tilde{y}$. Then there exists $t_{xy} \in (0,1)$ such that $\gamma(t)=\tilde{\gamma}(t)$ for any $t \in [0,t_{xy}]$, see Figure 2. 
\end{proof}

According to \cite[Definition 4.5]{Sturm},  a complete separable geodesic metric measure space $(X,d,\meas)$ satisfies the $\CD(K,\infty)$ condition,  where $K \in \mathbb{R}$, if for any $\mu_0, \mu_1 \in \mathcal{P}_2^{\text{ac}}(X,\meas)$ there exists a Wasserstein geodesic $(\mu_t)_{t \in [0,1]}$ between them such that
\begin{equation}\label{eq:conv}
\mathrm{Ent}(\mu_t) \le (1-t) \mathrm{Ent}(\mu_0) + t\mathrm{Ent}(\mu_1) - K t (1-t) W_2^2(\mu_0,\mu_1)/2
\end{equation}
for any $t \in [0,1]$, where the relative entropy  $\mathrm{Ent}$ is defined on $\cP(X)$ by
\[
 \mathrm{Ent} (\mu) =  \begin{cases}
 \int_X\ln \rho \di \mu & \text{ if } \mu = \rho \meas \in \mathcal{P}_2^{\text{ac}}(X,\meas),\\
 +\infty & \text{otherwise.}
 \end{cases}
\]

\begin{prop}
The space $(\cC,\widehat g,\widehat \cH^2)$ satisfies no $\CD(K,\infty)$ condition, whatever $K \in \mathbb{R}$.
\end{prop}

\begin{proof}
If $\mu_0,\mu_1$ are as in the proof of the previous proposition,  they both have finite relative entropy.  Moreover, any Wasserstein geodesic $\mu = (\mu_t)_{t \in [0,1]}$ between them is such that for any $t \in [t_\eps,1-t_\eps]$ the interpolating measure $\mu_t$ is supported in $\Gamma_t = \{\gamma(t) : \gamma \in \Geo(X) \text{ with } \gamma(0) \in A_0 \text{ and } \gamma(1) \in A_1\}$ which is $\widehat \cH^2$-negligible. As a consequence,  for any such a $t$, the measure $\mu_t$ has infinite relative entropy, hence \eqref{eq:conv} cannot hold.
\end{proof}

\subsection{The example is degenerate and not finite-dimensional MCP}

Let $(X,d)$ be a complete separable geodesic metric space endowed with a Radon measure $\meas$ which is finite on bounded sets. We  say that $(X,d,\meas)$ is non-degenerate if for any Borel set $A \subset X$ such that $\meas(A)>0$, it holds that $\meas(A_{x,t})>0$ for any $x \in X$ and $t\in [0,1]$, where
\[
A_{x,t} = \{\gamma(t) : \gamma \in  \mathrm{Geo}(X) \text{ with } \gamma(0) \in A \text{ and } \gamma(1) = x \}.
\]
The non-degeneracy condition was proposed in \cite[Definition 4.1]{Kell} after \cite[Definition 3.1]{CavallettiMondino}.  For brevity, we say that $(X,d,\meas)$ is degenerate if it is not non-degenerate.

\begin{prop}
The space $(\cC,\widehat g,\widehat \cH^2)$ is degenerate. 
\end{prop}
\begin{proof}
Set $x = (0,\pi/2)$.   For $\eps\in (0,0.1)$,  consider $o=(\eps,0) \in \Omega_+$ and the set $A = B_{\eps}(o) \subset \Omega_+$.  Then there exists $t_\eps \in (0,1)$ such that for any $\gamma \in \Geo(X)$ with $\gamma(0) \in A$ and $\gamma(1) = x$,  the point $\gamma(t)$ belongs to $\partial \Omega_+$ for any $t  \in (t_\eps,1]$. As a consequence, $\widehat \cH^2(A_{x,t})=0$ for any such a $t$, see Figure 3.
\end{proof}

Introduced independently in \cite{SturmII} and \cite{Ohta}, the measure contraction property $\MCP(K,N)$, where $K \in \mathbb{R}$ and $N \in [1,+\infty)$, is a weak version of the finite-dimensional $\CD(K,N)$ condition. For the sake of brevity, we do not provide the definitions of these two conditions, but we refer the reader to \cite[Section 6.2]{CavallettiMilman}, for instance, for a nice account on them.

\begin{cor}
The space $(\cC,\widehat g,\widehat \cH^2)$ satisfies no mesure contraction condition $\mathrm{MCP}(K,N)$, whatever $K \in \R$ and $N \in [1,+\infty)$.
\end{cor}

\begin{proof}
For $A_{x,t}$ as in the proof of the previous proposition, we get from \cite[Lemma 2.3]{Ohta} that any finite dimensional $\MCP$ condition implies a positive lower bound on $\widehat \cH^2(A_{x,t})$ that would contradict $\widehat \cH^2(A_{x,t})=0$ for any $t  \in (t_\eps,1]$.
\end{proof}

\begin{figure}
\includegraphics[height=3.5cm]{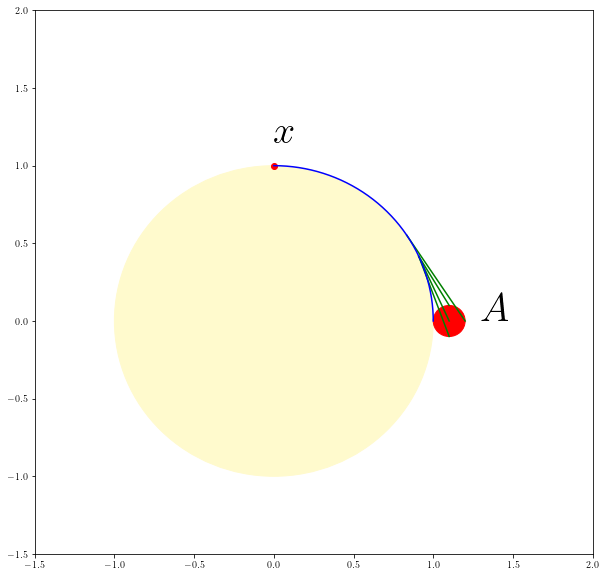}
\caption{Negligible set of $t$-intermediary points $\gamma(t)$ for $t$ close to $1$.}
\end{figure}

\section{Proof of Theorem \ref{prop:NotLp}}\label{sec:B}

This section is devoted to proving Theorem \ref{prop:NotLp}, which exhibits a compact example of a strong Kato limit which cannot be obtained as Gromov--Hausdorff limit of non-collapsing manifolds with a small $L^p$ bound on the Ricci curvature. 

Let $g_0$ be the round metric on  the sphere $\mathbb{S}^2$. In Fermi coordinates $(t,\theta) \in$ $(-\pi/2,\pi/2) \times \mathbb{S}^1$ around an equator $\mathcal{E} = \{t=0\}$, the round metric writes as
\begin{equation}\label{eq:round}
g_0 = \di t^2 + \cos^2 t \, \di \theta^2.
\end{equation}
For $\alpha \in (0,\pi/2)$, set
\[
\overline{g} := e^{2\overline{u}} g_0 \qquad \text{with} \quad \overline{u}(t,\theta) = \log\left(\frac{\cos (\alpha)}{1 - \sin (\alpha) \sin (|t|)}\right).
\]
See Figure 4. The following is a direct consequence of Theorem \ref{th:IBC}.

\begin{prop}
The space $(\mathbb{S}^2,\overline{g})$ is a non-collapsed strong Kato limit.
\end{prop}

\begin{proof}
Observe that $\overline{u}$ is continuous, bounded and $g_0$-Lipschitz.  Moreover, in Fermi coordinates, the Laplacian $\Delta_{g_0}f$ of a function $f$ only depending on $t$ writes as
$$\Delta_{g_0}f = -\partial_{tt}^2 f+\tan(t)\partial_t f.$$
A direct computation for $f=\overline{u}$ gives that, in the sense of distributions, 
\begin{equation}\label{eq:Laplacien}
\Delta_{g_0}\overline{u}= - 2  \sin (\alpha)  \, \di\mathcal{H}_{g_0}^1 \measrestr \mathcal{E} - \frac{-\sin^2(\alpha) +2\sin(\alpha)\sin(|t|)-\sin^2(\alpha)\sin^2(|t|)}{[1-\sin(\alpha) \,   \sin(|t|)]^2} \, \di v_{g_0}
\end{equation}
where $\di \mathcal{H}_{g_0}^1$ is the one-dimensional Hausdorff measure of $g_0$. As a consequence, $\Delta_{g_0}\overline{u}$ is a bounded Radon measure. Therefore, we can apply Theorem \ref{th:IBC} to conclude. 
\end{proof}

In the remainder of this section we show that $(\mathbb{S}^2,\overline{g})$ cannot be obtained as the limit of a sequence of non-collapsing surfaces $\{ (M_\ell, g_\ell)\}$ satisfying a uniform smallness $L^p$ bound in the sense of Petersen-Wei. For this, we start by recalling some notations, the relation between an $L^p$ and a strong Kato bound, and some results about surfaces with bounded integral curvature. 

\subsection{Strong Kato and $L^p$ bounds.} Let us briefly recall some notions from \cite{PW1,PW2} with the point of view adopted in \cite{Aubry}. For any $n \in \mathbb{N}$, $\kappa \ge 0$ and $r>0$,  set
\[
\mathbb{V}_{n,\kappa}(r) := \int_0^r \left(\frac{\sinh(\kappa t)}{\kappa}\right)                                                                                                                                                                                                                                                                                                                                                                                                                                                                                                        ^{n-1} \di t \,\,\, \text{ for $\kappa>0$}, \qquad \mathbb{V}_{n,0}(r) := \frac{r^n}{n} \, \cdot
\]
Note that $\mathbb{V}_{n,\kappa}(r)$ is a multiple of the Riemannian volume of a ball of radius $r$ in the $n$-dimensional space form $\mathbb{M}^n_{-\kappa^2}$ of constant sectional curvature equal to $-\kappa^2$. 

\begin{figure}
\includegraphics[height=3cm]{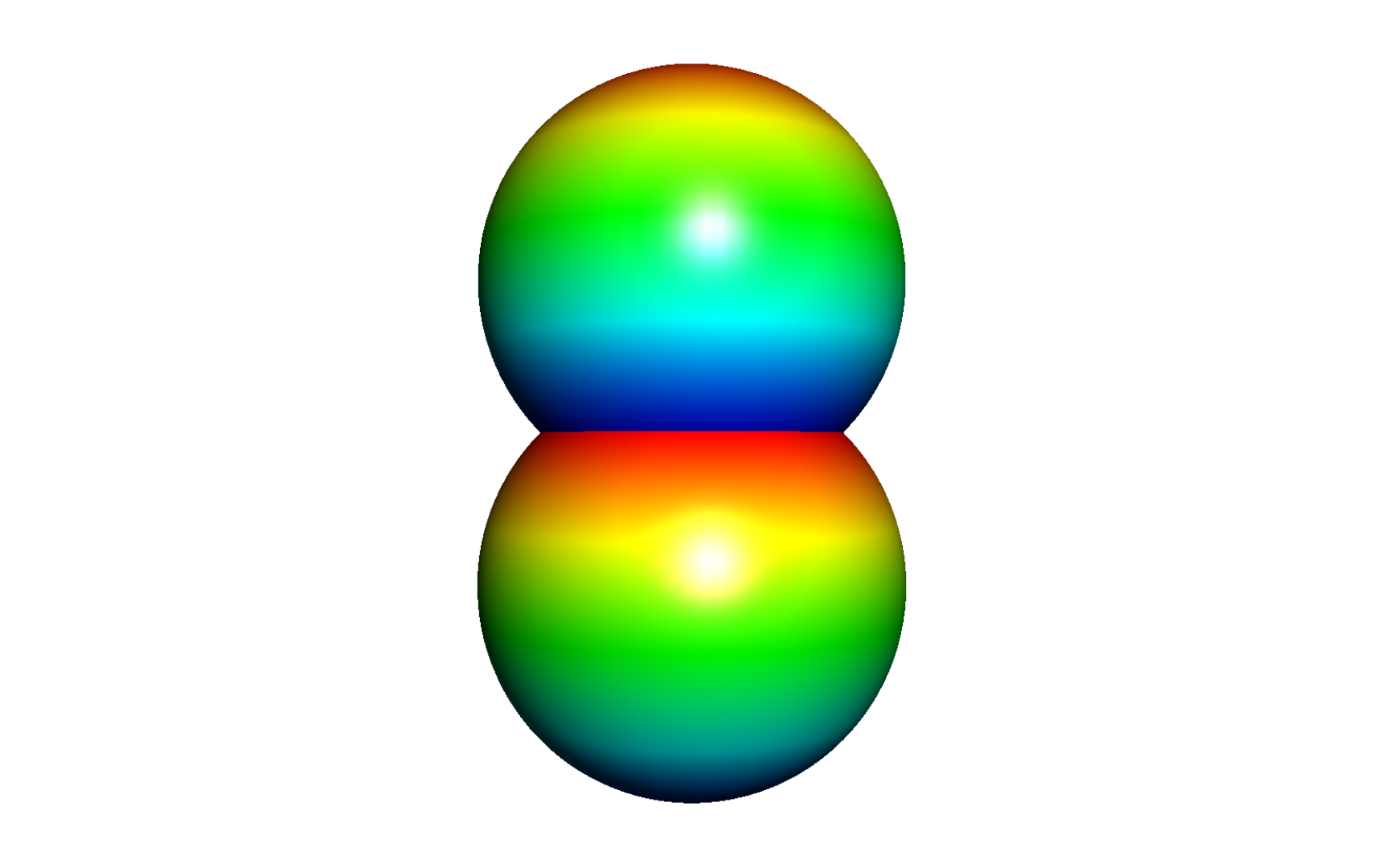}
\caption{The space $(\mathbb{S}^2,\overline{g})$.}
\end{figure}

On a closed smooth Riemannian manifold $(M^n,g)$ of diameter at most $D$, for fixed $\kappa \geq 0$ define the non-negative function
$$\uprho_\kappa=(\ricci+(n-1)\kappa^2g)_-$$
which represents the amount of Ricci curvature below $-(n-1)\kappa^2$.  For fixed $p > \frac{n}{2}$, a smallness $L^p$ bound means that $\uprho_\kappa$ is uniformly controlled in $L^p$ by a small constant,  i.e.
\begin{equation}\label{eq:assumptionLp}\tag{$L^p_{-\kappa^2}$}
D^2 \left( \fint_M \uprho_\kappa^p  \di\vol_g \right)^{\frac{1}{p}} \le \epsilon(n,p) \left( \frac{D^n}{\mathbb{V}_{n,\kappa}(D)}\right)^{\frac{1}{p}},
\end{equation}
where $\epsilon(n,p)$ is related to the constant $B(n,p)$ appearing in \cite[Théorème 4.6]{Aubry} by
\[
\left( \frac{1}{\epsilon(n,p)}\right)^{\frac{p}{2p-1}} = 4 B(p,n) 2^{\frac{1}{2p-1}}.
\]
We know from \cite[Theorème 4.8 and Lemme 4.10]{Aubry} that \eqref{eq:assumptionLp} implies
\begin{equation}\label{eq:prop1b}
\mathrm{v}_g(B_g(x,r)) \le C(n,p) \mathbb{V}_{n,\kappa}(r)
\end{equation}
for some $C(n,p)>0$ depending only on $n$ and $p$, and
\begin{equation}\label{eq:prop2b}
r^2 \left( \fint_{B(x,r)} \uprho_\kappa^p  \di\vol_g\right)^{\frac{1}{p}} \le 2 \left( \frac{r}{D} \right)^2 \left( \frac{\mathbb{V}_{n,\kappa}(D)}{\mathbb{V}_{n,\kappa}(r)} \right)^{\frac 1p}  D^2 \left( \fint_M \uprho_\kappa^p  \di\vol_g \right)^{\frac 1p}
\end{equation}
for any $x \in M$ and $r \in (0,D)$.  Applying \eqref{eq:assumptionLp} further to the previous yields
\begin{equation}\label{eq:further}
r^2 \left( \fint_{B(x,r)} \uprho_\kappa^p  \di\vol_g\right)^{\frac{1}{p}} \le 2 \epsilon(n,p) \left( \frac{r}{D} \right)^2 \left( \frac{D^n}{\mathbb{V}_{n,\kappa}(r)} \right)^{\frac 1p} \, \cdot
\end{equation}

It follows from \cite{CarronRose} that the assumption \eqref{eq:assumptionLp} implies a strong Kato bound (see also \cite{RoseStollmann} for a related statement). For the sake of completness, we provide a proof. 

\begin{lem}\label{lem:LpSK}
Let $(M^n,g)$ be a closed manifold of diameter at most $D$. Under \eqref{eq:assumptionLp}, there exist $T>0$ depending only on $n,p,\kappa,D$ and $C_0(n,p)>0$ depending only on $n, p$, such that 
\begin{equation}
\label{eq:Dynkin}k_T(M,g) \leq \frac{1}{16n}
 \end{equation}
and for any $t \in (0,T)$,
\begin{equation}
\label{eq:SK}
\tag{SK}
 k_t(M,g)\leq C_0(n,p)\left(\frac{t}{T}\right)^{1- \frac{n}{2p}}. 
 \end{equation}
  
\end{lem} 
\proof
In this proof, $C_0(n,p)$ denotes a positive constant depending on $n,p$ only whose value may change from line to line. From \eqref{eq:further}, for any $x \in M$ and $r\in (0,D)$,
\begin{align*}
r^2\fint_{B_g(x,r)} \Ric_{-}\di\vol_g &
\leq r^2(n-1)\kappa^2+r^2\left(\fint_{B_g(x,r)} \uprho_\kappa^p\di\vol_g \right)^{\frac 1p} \\
& \leq r^2(n-1)\kappa^2+ 2\epsilon(n,p)\left( \frac{r}{D} \right)^2 \left( \frac{D^n}{\mathbb{V}_{n,\kappa}(r)}\right)^{\frac 1p},  
\end{align*}
and since $\mathbb{V}_{n, \kappa}(r) \geq r^n/n$ we get
\begin{equation}
\label{eq:prop3}
r^2\fint_{B_g(x,r)} \Ric_{-}\di\vol_g \leq r^2(n-1)\kappa^2+C_0(n,p)\left( \frac{r}{D} \right)^{2-\frac{n}{p}}. 
\end{equation}
Theorem 4.3 in \cite{CarronRose} states that there exists $\updelta_n>0$ only depending on $n$ such that  if
$$\int_0^D e^{-\frac{r^2}{7R^2}}r \left(\fint_{B_g(x,r)}\Ric_-\di\vol_g \right) \di r \leq \updelta_n$$
for some $R \in (0,D]$, then $k_{R^2}(M,g) \leq \frac{1}{16n}$. Because of \eqref{eq:prop3} we have
\begin{equation}
\label{eq:prop4}
\int_0^D e^{-\frac{r^2}{7R^2}}r \left(\fint_{B_g(x,r)}\Ric_-\di\vol_g \right) \di r \leq 25(n-1)\kappa^2R^2+C_0(n,p)\left(\frac{R}{D} \right)^{2-\frac{n}{p}}.
\end{equation}
Set $\updelta(n,p)=\left(\frac{\updelta_n}{2C_0(n,p)}\right)^{1/(2-n/p)}$. Then by choosing 
$$R_0= \min \left\{ \frac{\sqrt{\updelta_n}}{\kappa\sqrt{50(n-1)}}, \updelta(n,p) D\right\},$$
and $T=R_0^2$ we obtain \eqref{eq:Dynkin}. Moreover, Corollary 4.1 of \cite{CarronRose} ensures that there exists $\uplambda$ only depending on $n$ and $p$ such that for any $R \leq D$
$$k_{R^2}(M,g) \leq \uplambda \int_0^D e^{-\frac{r^2}{7R^2}}r \left(\fint_{B_g(x,r)}\Ric_-\di\vol_g \right) \di r.$$
Therefore, inequality \eqref{eq:prop4} also implies that \eqref{eq:SK} holds for any $t< T$.
\endproof

\subsection{Convergence of surfaces with Bounded Integral Curvature.}
We shall use convergence results for closed surfaces due to C.~Debin and A.~D.~Alexandrov. We very briefly recall some notation to give the statement that we need. 
We refer to \cite{Debin, Troyanov} for the precise definition of metrics with Bounded Integral Curvature, BIC for short. On a closed surface $M$, a metric $\dist$ with BIC is associated to a curvature measure $\omega$ which, in the case of a smooth Riemannian metric $g$ on $M$, coincides with $\omega_g=K_g \di\vol_g$. The contractibility radius is the largest $r>0$ such that all closed balls of radius $s<r$ are homeomorphic to a closed disc: see \cite[Section 2]{Debin} for a description of its properties. In  the next statement, we combine the main theorem of \cite{Debin} and a result of Alexandrov \cite{Aleks67}, see also \cite[Theorem 2.6]{Troyanov}), for the convergence of the curvature measures. 

\begin{thm}
\label{thm:Debin}
Let $M$ be a closed surface and $\{\dist_\ell\}_{\ell\in \N}$ be a sequence of metrics with \emph{BIC}. Let $A, c, R, \delta$ be positive constants. Assume that for any $\ell$ :
\begin{enumerate}
\item for any $x \in M$, the positive part of the curvature measure $\omega_\ell$ associated to $\dist_\ell$ satisfies
$$\omega_\ell^+(B_{\dist_\ell}(x,R))\leq 2\pi-\delta;$$
\item the contractibility radius of $\dist_\ell$ is bounded below by $c$; 
\item the $2$-dimensional Hausdorff measure of $\dist_\ell$ satisfies $\mathcal{H}^2_{\ell}(M)\leq A$.
\end{enumerate}
Then there exists a metric $\dist$ with \emph{BIC} such that, up to passing to a subsequence, there are diffeomorphisms $\phi_\ell: M \to M$ such that $\phi_\ell^*\dist_\ell$ converges uniformly to $\dist$ on $M$. 

If moreover there exists $C>0$ such that $|\omega_\ell|(M)\leq C$ for all $\ell \in \N$,  then the curvature measure $\omega$ of $\dist$ is the weak limit of $\omega_\ell$, that is for any $f \in C^0(M)$ we have
$$\lim_{\ell \to \infty}\int_{M} f \di \omega_\ell= \int_M f\di \omega. $$
\end{thm}

In the case of $(\mathbb{S}^2, \overline{g})$, the conformal transformation rule  for the Gauss curvature \eqref{eq:CT} and the expression of $\Delta_{g_0}\overline{u}$ given in  \eqref{eq:Laplacien} ensure that the curvature measure of $\overline{g}$ writes as
\begin{equation}
\label{eq:curvmeas}
\omega_{\overline{g}} = - 2  \sin (\alpha) \mathcal{H}^1_{\overline{g}} \, \measrestr \mathcal{E}  + \mathrm{v}_{\overline{g}},
\end{equation}
where $ \cH_{\overline{g}}$ is the one-dimensional Hausdorff measure of $\overline{g}$. In particular,  the negative part of $\omega_{\overline{g}}$ is singular and concentrated in the equator:
\[
(\omega_{\overline{g}})^{-} = 2  \sin (\alpha) \mathcal{H}^1_{\overline{g}} \,\measrestr \mathcal{E}.
\]

\subsection{Proof of Theorem \ref{prop:NotLp}}

We need the following lemma, which mainly relies on the volume bound \eqref{eq:prop1b} and on a formula, first given in the work of F.~Fiala, for the variation of the length of a geodesic disk on a surface.

\begin{lem}\label{lem:curvaturebound}
Let $(M^2,g)$ be a closed surface of diameter at most $D$ such that \eqref{eq:assumptionLp} holds for some $p>1$ and $\kappa \ge 0$.  Consider the measures $\omega_g^{\pm} := (K_g)_{\pm} \mathrm{v}_g$. Then there exists $R_1>0$ depending on $p,\kappa,D$ only such that for any $x \in M$ and $r \in (0,R_1)$,
\[
\omega_g^+(B_g(x,r/2)) \le 8\pi\left(1 - \frac{\vol_g(B_g(x,r))}{\pi r^2}\right) + \frac{\pi}{100} \, \cdot 
\]
\end{lem}

\begin{proof}
Consider $x \in M$ and $R \in (0,D)$.  Using \cite[Theorem 3, Theorem 5]{Fiala}, see also \cite[Theorem 6.1]{Hartman} and \cite{ST}, with Gauss-Bonnet's formula, we get
\begin{equation}
\label{eq:Fiala}
\vol_g(B_g(x,R)) \leq  2\pi \int_0^R\chi(B_g(x,r))(R-r)\di r- \int_0^R (R-r) \left(\int_{B_g(x,r)} K_g \di\vol_g \right) \di r,
\end{equation}
where  $\chi(B_g(x,r))$ is the Euler characteristic of the ball $B_g(x,r)$. This implies
\begin{equation}
\label{eq:Fiala1}
\vol_g(B_g(x,R)) \leq  \pi R^2- \int_0^R (R-r) \left(\int_{B_g(x,r)} K_g \di\vol_g \right) \di r
\end{equation}
because $\chi(B_g(x,r))\leq 1$.  Since
\begin{equation*}
\int_0^R (R-r) \left(\int_{B_g(x,r)} (K_g)_- \di\vol_g \right) \di r \le \omega_g^-(B_g(x,R) )\int_0^R (R-r) \di r  = \frac{R^2}{2} \omega_g^-(B_g(x,R))
\end{equation*}
and
\begin{align*}
\int_{R/2}^R (R-r) \left(\int_{B_g(x,r)} (K_g)_+ \di\vol_g \right) \di r \ge  \omega_g^+(B_g(x,R/2) )\int_{R/2}^R (R-r) \di r  = \frac{R^2}{8} \omega_g^+(B_g(x,R/2))
\end{align*}
the previous yields that
\[
\frac{R^2}{8} \omega_g^+(B_g(x, R/2)) \le \pi R^2 - \vol_g(B_g(x,R)) + \frac{R^2}{2} \omega_g^-(B_g(x,R)).
\]
Now we use Hölder's inequality,  \eqref{eq:further} and \eqref{eq:prop1b},  then $\mathbb{V}_{2,\kappa} (R) \leq \frac{R^2}{2}\cosh(\kappa R)$, to obtain:

\begin{align*}
\omega_g^-(B_g(x,R)) & \leq \int_{B_g(x,R)} (K_g+\kappa^2)_-\di\vol_g + \kappa^2\vol_g(B_g(x,R)) \\
&  \leq\vol_g(B_g(x,R)) \left( \fint_{B_g(x,R)} (K_g+\kappa^2)_-^p \di \vol_g\right)^{\frac{1}{p}}+ \kappa^2\vol_g(B(x,R)) \\
& \leq 2 C(2, p)\, \epsilon(2,p)\mathbb{V}_{2,\kappa} (R)^{1-\frac{1}{p}} D^{\frac{2}{p}-2}+ C(2,p)\kappa^2\mathbb{V}_{2,\kappa}(R) \\
& \leq 2^{\frac{1}{p}}C(2,p)\, \epsilon(2,p) \left(\frac{R}{D}\right)^{2-\frac{2}{p}}\cosh(\kappa R)^{1-\frac 1p}+ \frac{C(2,p)\kappa^2 R^2}{2}\cosh(\kappa R).
\end{align*}
Therefore we get
\begin{align*}
\frac{R^2}{8} \omega_g^+(B_g(x,R/2))&  \leq (\pi R^2 - \vol_g(B_g(x,R))) + 2^{\frac{1}{p}-2}C(2,p)\, \epsilon(2,p) R^2 \left(\frac{R}{D}\right)^{2-\frac{2}{p}}\cosh(\kappa R)^{1-\frac 1p}\\
& + \frac{C(2,p)\kappa^2 R^4}{4}\cosh(\kappa R).
\end{align*}
By choosing $R_1$ small enough depending on $p, \kappa, D$ we obtain the desired inequality. \end{proof}

\begin{rem}\label{rem:borne_neg_curv_meas}
The previous proof also implies the following result. Consider a closed surface $(M^2,g)$ of diameter at most $D$ satisfying \eqref{eq:assumptionLp} for some $p>1$ and $\kappa \ge 0$.  Then there exists $\Lambda(p,\kappa,D)>0$ depending only on $p,\kappa,D$ such that for any $x \in M$ and $r \in (0,D)$,
\begin{equation}
\label{eq:lemmaB}
\omega_g^{-}(B_g(x,r)) \le \Lambda(p,\kappa,D)r^{2-\frac{2}{p}}.
\end{equation}

\end{rem}

We are now in position to prove Theorem \ref{prop:NotLp}.

\begin{proof}[Proof of Theorem \ref{prop:NotLp}]
Assume by contradiction that there exists a non-collapsing sequence $\{(M_\ell, g_\ell)\}_\ell$ of closed surfaces with diameters at most $D>0$, satisfying \eqref{eq:assumptionLp} for uniform $p>1$ and $\kappa\geq 0$,  such that $\{(M_\ell, g_\ell)\}_\ell$ converges to $(\mathbb{S}^2, \overline{g})$ in the Gromov--Hausdorff topology. Througout the proof,  we will use that Lemma \ref{lem:LpSK} implies that this sequence also satisfies the uniform strong Kato bound \eqref{eq:SK} . 

\textbf{Step 1.} We start by proving topological stability for the sequence, i.e.~$M_\ell$ is homeomorphic to $\mathbb{S}^2$ for any large enough $\ell$. 

From \cite{CMT23}, there exist a conformal metric $\tilde{g}_\ell$ on each $M_\ell$ such that the Gauss curvature of $\tilde{g}_\ell$ is uniformly bounded from below by some negative constant $C$ depending only on $p,\kappa,D$.  As a consequence,  there exists a distance $\tilde{\dist}$ on $\mathbb{S}^2$ such that $(\mathbb{S}^2,\tilde{\dist})$ is an Alexandrov space to which $\{(M_\ell, \tilde{g}_\ell)\}$ converges in Gromov-Hausdorff topology, up to extraction of a subsequence. Then Perelman's topological stability for Alexandrov spaces (see \cite{Kapovitch}) implies that $M_\ell$ is a topological sphere for $\ell$ large enough. 

From now on, we only consider $\ell$ such that $M_\ell=\mathbb{S}^2$. 

\textbf{Step 2.} We show that there is $C>0$ only depending on $p,\kappa,D$ such that for large $\ell$,
\begin{equation}\label{eq:step2}
\int_{\mathbb{S}^2}|K_{g_\ell}|\di\vol_{g_\ell} \leq C.
\end{equation}

The Gauss-Bonnet theorem, Hölder's inequality,  and \eqref{eq:further}, successively imply that
\begin{align*}
\int_{\mathbb{S}^2} (K_{g_\ell})_{+} \di\vol_{g_\ell} & = \int_{\mathbb{S}^2} (K_{g_\ell})_{-} \di\vol_{g_\ell} +4\pi \\
& \leq \vol_{g_\ell}(\mathbb{S}^2)^{\frac 1q} \left(\int_{\mathbb{S}^2}  (K_{g_\ell})^p_{-} \di\vol_{g_\ell} \right)^{\frac 1p}+4\pi \qquad \text{with $q = p/(p-1)$}\\
&\leq \vol_{g_\ell}(\mathbb{S}^2)  \left(\fint_{\mathbb{S}^2}  (K_{g_\ell}+\kappa^2)^p_{-} \di\vol_{g_\ell} \right)^{\frac 1p} + \kappa^2 \vol_{g_\ell}(\mathbb{S}^2) +4\pi\\
&\leq \epsilon(2,p)C(2,p)D^{\frac{2}{p}} \mathbb{V}_{2,\kappa}(D)^{1-\frac{1}{p}} + C(2,p)\kappa^2 \mathbb{V}_{2,\kappa}(D) +4\pi.
\end{align*}
We eventually get existence of $C>0$ depending on $p,\kappa,D$ such that \eqref{eq:step2} holds.

\textbf{Step 3.} We show that the assumptions of Theorem \ref{thm:Debin} are satisfied. 

By the previous step, we already know that the curvature measures $\omega_\ell=K_{g_\ell}\di\vol_{g_\ell}$ are uniformly bounded.  By continuity of the conformal factor $e^{2\overline{u}}$ in the definition of $\overline g$, there exists $R_2>0$ such that for any $r \in (0,R_2)$ and $x \in \mathbb{S}^2$,
\[
\vol_{\overline{g}} (B_{\overline{g}}(x,r)) \ge \pi r^2 - \frac{\pi}{200} r^2.
\]
Let $R_1$ be as in Lemma \ref{lem:curvaturebound}, and  $R \in (0, \min\{R_1,R_2\})$. By volume continuity \cite[Section 7]{CMT} there exists $\ell_0$ depending on $R$ such that whenever $\ell \ge \ell_0$, for any $x \in \mathbb{S}^2$, 
\begin{equation}\label{eq:bound}
\vol_{g_\ell} (B_{g_\ell}(x,R)) \ge \pi R^2 - \frac{\pi}{100} R^2.
\end{equation}
Then Lemma \ref{lem:curvaturebound} implies that for any $\ell \ge \ell_0$ and $x \in M_\ell$,
\begin{equation}
\label{eq:omega+}
\omega_{g_\ell}^+(B_{g_\ell}(x,R)) \le \frac{8\pi}{100} + \frac{\pi}{100} \le \frac{\pi}{10}, 
\end{equation}
thus the first assumption of Theorem \ref{thm:Debin} holds. 

As for the second assumption, we aim to prove that there exists a uniform positive lower bound on the contractibility radius of $(\mathbb{S}^2, g_\ell)$.  In the following, we omit the dependence on $p, \kappa, D$ of $\Lambda$ appearing in equality \eqref{eq:lemmaB}.  Choose $R>0$ small enough to ensure that
\[
R < \min\{R_1,R_2\}\quad \text{and} \quad \frac{99}{100} \pi -\frac{p^2 \Lambda}{\pi} R^{2-\frac{2}{p}} > \frac{1}{2} \, \cdot
\]
Using inequalities \eqref{eq:Fiala} and \eqref{eq:lemmaB} we obtain that for any $\ell \geq \ell_0$
\begin{align*}
\vol_{g_\ell}(B_{g_\ell}(x,R)) & \leq 2\pi \int_0^{R}\chi(B_{g_\ell}(x,r))(R-r)\di r+\int_0^{R}(R-r)\left( \int_{B_{g_\ell}(x,r)}(K_{g_\ell})_{-} \di \vol_{g_\ell} \right) \di r \\
& \leq 2\pi \int_0^{R}\chi(B_{g_\ell}(x,r))(R-r)\di r+p^2 \Lambda R^{4-\frac{2}{p}}.
\end{align*}
Since $R < \min\{R_1,R_2\}$ we use \eqref{eq:bound} to get that for any $\ell \ge \ell_0$,
\[
\frac{99}{100}\pi R^2 \le 2\pi \int_0^{R}\chi(B_{g_\ell}(x,r))(R-r)\di r+p^2 \Lambda R^{4-\frac{2}{p}} \,\cdot 
\]
Set $r_\ell(x) := \sup \{r>0 : \chi(B_{g_\ell}(x,s)) = 1 \text{ for a.e.~}s \in (0,r)\}$.  Observe that, by definition of $r_\ell(x)$, there is a set of null measure of radii $r \in (0,r_\ell(x))$ for which $B_{g_\ell}(x,r)$ might not be homeomorphic to a Euclidean disk. For any such radius $r$, there exists $r_0 \in (r, r_\ell(x))$ so that $B_{g_\ell}(x,r_0)$ has Euler characteristic equal to 1, thus it is homeomorphic to a disk. This together with inequality \eqref{eq:omega+} ensures that we can apply \cite[Lemma 2.3]{Debin} to get that for any $s \in (0,\min\{r_0,R\}]$ the ball $B_{g_\ell}(x,s)$ is homeomorphic to a disk. Then for all $r \in (0, \min\{r_\ell(x),R\})$ the ball $B_{g_\ell}(x,r)$ is homeomorphic to a disk. Therefore, a lower bound on $r_\ell(x)$ will give a lower bound on the contractibility radius. Moreoever,
$\chi(B_{g_\ell}(x,r) )\le 0$ for a.e.~$r>r_\ell(x)$. If $r_\ell(x) < R$, we get
\begin{align*}
\frac{99}{100} \pi R^2 & \le 2\pi \int_0^{r_\ell(x)}(R-r)\di r+p^2 \Lambda R^{4-\frac{2}{p}}\le 2\pi R r_\ell(x)+p^2 \Lambda R^{4-\frac{2}{p}}
\end{align*}
so that
\[
r_\ell(x) \ge R \left( \frac{99}{100} \pi -\frac{p^2 \Lambda}{\pi} R^{2-\frac{2}{p}} \right) > \frac{R}{2} \, \cdot
\]
Otherwise $r_\ell(x) \ge R>R/2$.  In the end,  we get that the contractibility radius of $(\mathbb{S}^2, g_\ell)$ is uniformly bounded from below by $R/2$.

Lastly, it follows from \eqref{eq:prop1b} that for any $\ell$,
\[
\vol_{g_\ell}(M_\ell) \le C(n,p,\kappa,D).
\]

As a consequence, Theorem \ref{thm:Debin} ensures the existence of a BIC metric $\dist$ on $\mathbb{S}^2$ and of diffeomorphisms \mbox{$\phi_\ell : \mathbb{S}^2 \to \mathbb{S}^2$} such that, up to a subsequence, $(\phi_\ell)^*\dist_{g_\ell} \to \dist$ uniformly on $\mathbb{S}^2$ and $\omega_{(\phi_\ell)^*\dist_{g_\ell} } \weakto \omega_{\dist}$. Since the Gromov-Hausdorff limit of a sequence is unique up to isometry, there exists a homeomorphism $f$ such that $f^*\dist=\dist_{\overline g}$.  Therefore,  the distances $\tilde{d}_\ell = f^* \dist_{(\phi_\ell)^*g_\ell}$ satisfy
\[
\tilde{d}_\ell \to \dist_{\overline{g}} \quad \text{uniformly on $\mathbb{S}^2$,} \qquad  \omega_{\tilde{d}_\ell} \weakto \omega_{\overline{g}}.
\]

\textbf{Step 4.} We obtain a contradiction by using weak convergence of the curvature measures and the fact that $\omega_{\overline{g}}^{-}$ is singular with respect to the Riemannian volume. For $\varepsilon \in (0,\pi/2)$, consider tubular neighbourhoods $U^\varepsilon = \{|t|<\epsilon\}$ and $U^{2\varepsilon}=\{|t|<2\epsilon\}$ of the equator $\mathcal{E} = \{t=0\}$. Let $f_\epsilon : \mathbb{S}^2 \to [0,1]$ be a continuous function equal to 1 on $U^\varepsilon$ and vanishing outside $U^{2\varepsilon}$.  Thanks to the explicit expression \eqref{eq:curvmeas} for $\omega_{\overline{g}}$,  the dominated convergence theorem implies that
\[
\lim_{\varepsilon \to 0} \left(  - \int_{\mathbb{S}^2} f_\epsilon \di \omega_{\overline{g}} \right) = 4 \pi e^{2c_\alpha} \sin(\alpha)  =: C_0 > 0.
\]
Let us consider $\epsilon>0$ small enough to ensure that
\[
- \int_{\mathbb{S}^2} f_\epsilon \di \omega_{\overline{g}}  \ge \frac{C_0}{2} \,  \cdot 
\]
By weak convergence $\omega_{\tilde{d}_\ell} \weakto \omega_{\overline{g}}$ there exists $\ell_0$ such that for any $\ell \ge \ell_0$,
\begin{equation}\label{eq:contradict_1}
 \frac{C_0}{4} \le - \int_{\mathbb{S}^2} f_\epsilon \di \omega_{\tilde{d}_\ell}  \le   \int_{\mathbb{S}^2} f_\epsilon \di [\omega_{\tilde{d}_\ell}]^{-} \le [\omega_{\tilde{d}_\ell}]^{-}(U^{2\epsilon}) \,  \cdot 
\end{equation}
By triangle inequality there exist $N(\epsilon)$ points $\{x_i\}$ on $\mathcal{E}$ such that
$$
U^{2\epsilon} \subset \bigcup_{i} B_{g_0}(x_i,3\epsilon)
$$
and $N(\epsilon) \le \Gamma/\epsilon$ for some $\Gamma>0$ independent of $\epsilon$.  Since $g_0$ and $\overline{g}$ are bi-Lipschitz equivalent, and $\tilde{d}_\ell \to \dist_{\overline{g}}$ uniformly on $\mathbb{S}^2$, there exists $C \ge 3$ such that for any $\ell \ge \ell_0$,
$$
U^{2\epsilon} \subset \bigcup_{i} B_{\tilde{d}_\ell}(x_i,C\epsilon).
$$
By Remark \ref{rem:borne_neg_curv_meas}, we get that for any $\ell \ge \ell_0$,
\begin{equation}\label{eq:contradict2}
[\omega_{\tilde{d}_\ell}]^{-}(U^{2\epsilon}) \le \sum_i [\omega_{\tilde{d}_\ell}]^{-}(B_{\tilde{d}_\ell}(x_i,C\epsilon)) \le N(\epsilon)C(p,\kappa)\epsilon^2(1+\epsilon^{2-\frac{1}{p}}) \le \tilde{C}(p,\kappa) \epsilon
\end{equation}
for some $\tilde{C}(p,\kappa)>0$. We get a contradiction between \eqref{eq:contradict_1} and \eqref{eq:contradict2} by choosing a possibly smaller $\epsilon$ such that $ \tilde{C}(p,\kappa) \epsilon < C_0/4$. \end{proof}

\bibliographystyle{plain} 
\bibliography{BiblioBranch}
\end{document}